\documentclass[11pt]{article}
\usepackage{epsf,amssymb,latexsym,amsmath}

\newcommand{\proof}{\noindent{\bf Proof.\ }}
\newcommand{\qed}{\hfill $\square$ \bigskip}

\newtheorem{theorem}{\bf Theorem}[section]
\newtheorem{corollary}[theorem]{\bf Corollary}
\newtheorem{lemma}[theorem]{\bf Lemma}

\newtheorem{remark}[theorem]{\bf Remark}

\newtheorem{observation}[theorem]{\bf Observation}

\def\slika #1{\begin{center}\hskip 0.2mm\epsffile{#1}\end{center}}

\begin{document}

\title{A characterization of the edge connectivity of direct products of graphs}

\author{
Simon \v Spacapan 
\footnote{ University of Maribor, FME, Smetanova 17,
2000 Maribor, Slovenia. e-mail: simon.spacapan@uni-mb.si}
\footnote{This work is supported by Ministry of Education and Science of Slovenia, grants P1-0297 and J1-4106.}}
\date{\today}

\maketitle

\begin{abstract}
\noindent The direct product of graphs $G=(V(G),E(G))$ and $H=(V(H),E(H))$ is the graph, denoted as $G\times H$, 
with vertex set $V(G\times H)=V(G)\times V(H)$, where vertices $(x_1,y_1)$ and $(x_2,y_2)$ are adjacent in 
$G\times H$ if $x_1x_2\in E(G)$ and $y_1y_2\in E(H)$. The edge connectivity of a graph $G$, 
denoted as $\lambda(G)$, is the size of a minimum edge-cut in $G$. We introduce a function $\psi$ and 
prove the following formula 
$$\lambda (G\times H)=\min \{2\lambda(G)|E(H)|,2\lambda(H)|E(G)|,\delta(G\times H), 
\psi(G,H), \psi(H,G)\}\,.$$
We also describe the structure of every minimum edge-cut in $G\times H$.

\end{abstract}




\section{Introduction}
\begin{figure}[htb!]
\epsfxsize=14.3truecm
\slika{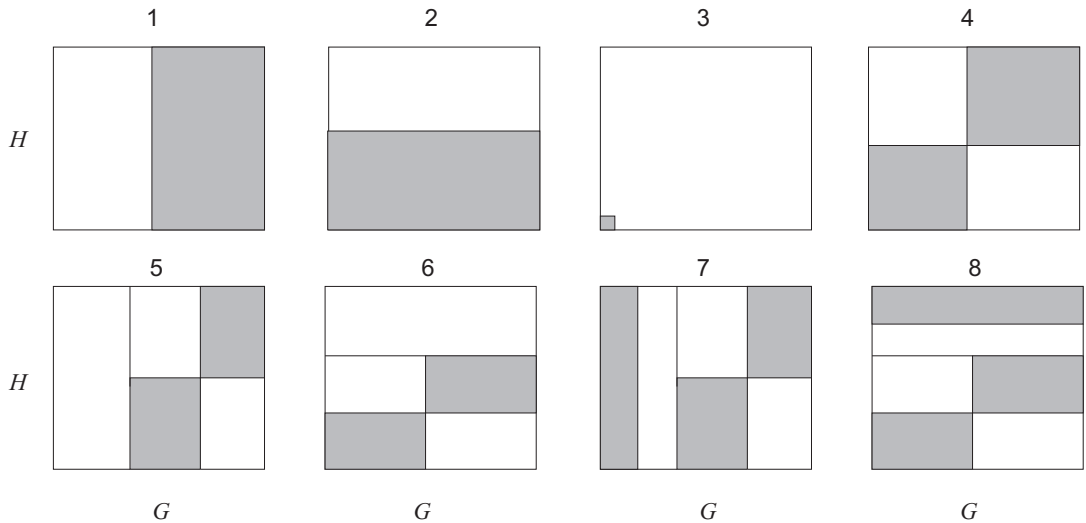}
\caption{Connected components of a minimum edge cut in $G\times H$.}
\label{struktura}
\end{figure}

Weichsel observed in \cite{we} that the direct product of graphs $G$ and $H$ is 
connected if and only if both graphs are connected and at least one of them is nonbipartite. 
Several decades after this observation a detailed study of connectivity of graph products 
was initiated by different authors. The aim of this study is to determine the connectivity of 
the product and express it in terms of connectivities of the factors. The second aim 
is to describe the structure and other properties of a minimum vertex-cut or a minimum edge-cut 
in the product. The size and the structure of a minimum edge-cut in the Cartesian product of graphs
was determined in \cite{jaz4} where the following result is proved 
$$\lambda (G\Box H)=\min \{\lambda(G)|V(H)|,\lambda(H)|V(G)|,\delta(G\Box H)\}\,.$$
The authors also prove that every minimum edge-cut in the Cartesian product of graphs is the preimage (under projection) 
of a minimum edge-cut of a factor, or the set of edges incident to a vertex of minimum degree. 
Therefore the connected components of a minimum edge-cut in the Cartesian product coincide 
with one of the first three cases 
of Fig. \ref{struktura}.
For the strong product of graphs a similar formula for the edge connectivity 
and an analogous result about the structure of a minimum edge-cut was proved, see \cite{jaz5}. 
A partial answer to these questions was obtained for the vertex connectivity of the 
Cartesian and the strong product of graphs, see \cite{jaz1,jaz2}. 

The problem of determining the connectivity of dirct products of graphs appears to be more challenging 
(compared to other graph products). 
This comes from the fact that the direct product of two bipartite graphs is not connected although 
both graphs may have high edge and vertex connectivities. Therefore it is not possible to 
express the connectivity of the direct product of graphs exclusively in terms of connectivities (and minimum degree, size and order)
of the factors. There is an extensive list of articles on connectivity of direct 
products of graphs where the authors consider special examples and determine their connectivity, 
or they obtain upper and lower bounds for the connectivity, or they study related concepts such as super connectivity, see 
\cite{jaz3, mi, rg, lg, lm1, jo, wy, ww1, ww, knjiga}. In this article we settle the question of edge connectivity 
of direct products as well as the question about the structure of a minimum edge-cut in the direct product of graphs. 

Before we state our main result we give the definitions and  the notation we use. 
Let $G$ and $H$ be graphs and $G\times H$ their direct product. 
For $x\in V(G)$ and $y\in V(H)$, the {\em $H$-layer} $H_x$, and the {\em $G$-layer} $G_y$, are defined as 
$$H_x=\{(x,h)\,|\,h\in V(H)\} ~{\rm and }~ G_y=\{(g,y)\,|\,g\in V(G)\}\,.$$

If $S$ is a subset of $E(G)$ then we call an edge of $S$ an {\em $S$-edge}. 
The graph obtained from $G$ by deleting all $S$-edges is denoted by $G-S$. 
A set $S\subseteq E(G)$ is an {\em edge-cut} in $G$ if $G-S$ is not connected, and 
the {\em edge connectivity} of $G$ is the size of a minimum edge-cut in $G$. 
If $S$ is a minimum edge-cut in $G\times H$ then we denote by $B$ and $W$ the connected components of $(G\times H)-S$.
We refer to $B$ as {\em black} and $W$ as {\em white}. 
For example, when we say that a vertex $(x,y)$ is white, we mean that $(x,y)\in W$. 
We say that an edge cut $S$ in $G\times H$ is of {\em type 1} if $B=B'\times V(H)$ and 
$W=W'\times V(H)$ where $B'\cup W'=V(G)$. 
The definitions of all other types of edge-cuts are evident from Fig. \ref{struktura}. 
For example $S$ is of {\em type 3} if $|B|=1$ or $|W|=1$. 
The {\em  bipartite edge frustration}  of a graph $G$, denoted as $\varphi (G)$, is the 
minimum number of edges whose delition makes $G$ a bipartite graph. The concept was introduced  
in \cite{bip} and studied in \cite{bip1,bip2,bip3}.
In section \ref{size} we formally define the function $\psi$ which is needed to determine the edge connectivity 
of direct products, moreover the bipartite edge frustration of factors is crucial in this definition. It follows from the 
definition that for a bipartite graph $H$, $\psi(G,H)$ is equal to the size of a minimum edge cut of type 4 or 5, 
and if $G$ is bipartite, then $\psi(H,G)$ is equal to the size of a minimum edge cut of type 4 or 6. 
We prove the following theorem in Section \ref{size}. 

\begin{theorem} \label{formula}
For any  graphs $G$ and $H$
$$\lambda (G\times H)=\min \{2\lambda(G)|E(H)|,2\lambda(H)|E(G)|,\delta(G\times H), 
\psi(G,H), \psi(H,G)\}\,.$$
\end{theorem}
The structure theorem is proved in Section \ref{zadnje}.

\begin{figure}[htb!]
\epsfxsize=10.5truecm
\slika{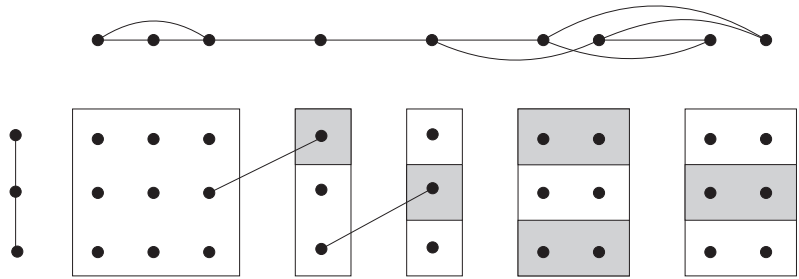}
\caption{An example of an edge-cut with no structure. }
\label{sporadic-example}
\end{figure}

\begin{theorem} \label{tipi}
If $G$ and $H$ are not equal to a path on three vertices or a 4-cycle then every minimum edge-cut 
in $G\times H$ is of type $1,\ldots,7$ or $8$.
\end{theorem}
It also follows from our results, that for every graph product $G\times H$ there is minimum 
edge-cut in $G\times H$ of type 1,2,3,4,5, or 6. In particulur, if there is a minimum edge-cut 
in $G\times H$ of type 8 (resp. 7), then there is also a minimum edge-cut of type  6 or 3 (resp. 5 or 3).
In Fig. \ref{sporadic-example} we give an example of a graph product that has a minimum edge-cut (of size 2)
with the structure that cannot be categorized as one of the eight types.

\section{The size of a minimum edge cut}
\label{size}

For each vertex $(x,y)\in V(G\times H)$ we denote by 
$K(x,y)$ the subgraph of $G\times H$ induced by 
$$(N_G(x)\times \{y\})\cup(\{x\}\times N_H(y))\,.$$
Note that  this is a complete bipartite subgraph of $G\times H$.

\begin{figure}[htb!]
\epsfxsize=5.0truecm
\slika{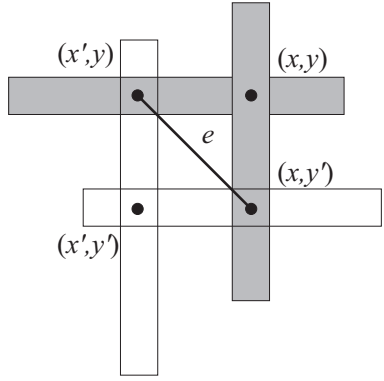}
\caption{$K(x,y)$ and $K(x',y')$ have a common edge $e=(x',y)(x,y')$.}
\label{lili}
\end{figure}

\begin{observation}\label{dva}
The edge $e=(x',y)(x,y')$ is an edge of $K(u,v)$ if and only if 
$(u,v)=(x,y)$ or $(u,v)=(x',y')$. 
\end{observation}
Observation \ref{dva} is depicted in Fig. \ref{lili}. The following observation is due to Weichsel \cite{we}.

\begin{observation} \label{kuku}
Let $G$ be a connected graph. If $G\times K_2$ is not connected then $G$ is 
bipartite and  connected components of $G\times K_2$ are 
$$C_1=(A\times \{x\})\cup (B\times \{y\})~{\text and }~C_2=(A\times \{y\})\cup (B\times \{x\})\,,$$
where $V(G)=A\cup B$ is the bipartition of graph G, and $V(K_2)=\{x,y\}$.
\end{observation}

\begin{lemma} \label{glavna}
Let $G$ and $H$ be graphs such that $\delta(G)\leq \delta (H)$. Let $S$ be a minimum edge cut in 
$G\times H$ and $B,W$ be connected components of $(G\times H)-S$. 
If $|S|<\delta(G\times H)$ then for every $(x,y)\in V(G\times H)$ either 
$N(x)\times \{y\}\subseteq B$ or $N(x)\times \{y\}\subseteq W$.

\end{lemma} 

\begin{figure}[htb!]
\epsfxsize=9.4truecm
\slika{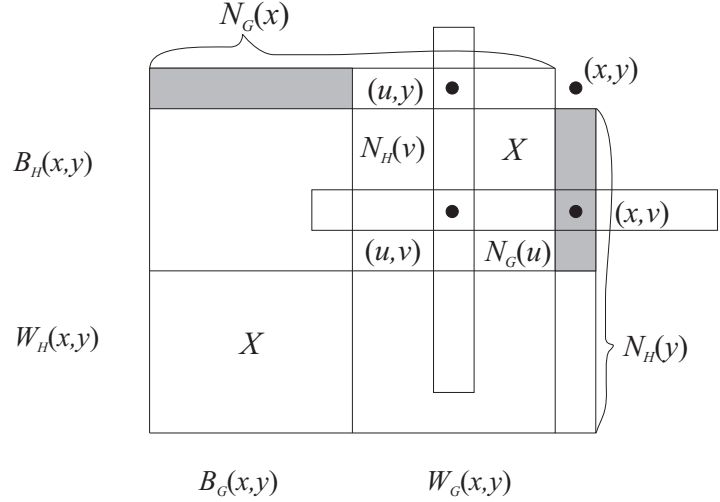}
\caption{The set $X$ and corresponding subgraphs $K(u,v),(u,v)\in X$.}
\label{glavna_slika}
\end{figure}

\proof Let $S$ be a minimum edge cut in $G\times H$ and let $B,W$ be
connected components of $(G\times H)-S$. 
We refer to $B$ as black and $W$ as white. Let $(x,y)$ be any vertex of $G\times H$ and 
let $B_G(x,y)$ and $W_G(x,y)$ be subsets of $N_G(x)$ such that 
$B_G(x,y)\times \{y\}$ is black and $W_G(x,y)\times \{y\}$ is white. Additionaly let 
$B_H(x,y)$ and $W_H(x,y)$ be subsets of $N_H(y)$ such that 
$\{x\}\times B_H(x,y)$ is black and $\{x\}\times W_H(x,y)$ is white (see Fig.~\ref{glavna_slika}).

We first prove the lemma under assumption that $\delta(H)\geq 3$.
Assume (contrary to the claim of the lemma) that $|B_G(x,y)|,|W_G(x,y)|\geq 1$. To prove the lemma we will show that 
$|S|\geq \delta(G\times H)$. To prove this we count the number 
of $S$-edges in each subgraph $K(u,v)$, where 
$$(u,v)\in (B_G(x,y)\times W_H(x,y))\cup (W_G(x,y)\times B_H(x,y))=X\,.$$

For each $(u,v)\in X$ let us denote by $\alpha'(u,v)$ the number of $S$-edges of $K(u,v)$ that have both endvertices in 
$B_G(x,y)\times W_H(x,y)$, or both endvertices in $W_G(x,y)\times B_H(x,y)$, or one endvertex in 
$W_G(x,y)\times W_H(x,y)$ and the other in $B_G(x,y)\times B_H(x,y)$. 
Note that $\alpha'(u,v)$ is the number of $S$-edges of $K(u,v)$ that are 
counted twice when counting $S$-edges of $K(u,v)$'s and 
$(u,v)$ goes through the set $X$ (see Observation \ref{dva}).
Let $\alpha''(u,v)$ be the number of all other $S$-edges of $K(u,v)$ and set
$$\alpha(u,v)=\frac12 \alpha'(u,v)+\alpha''(u,v)\,.$$
Clearly $$|S|\geq \sum_{(u,v)\in X}\alpha(u,v)\,.$$

\begin{figure}[htb!]
\epsfxsize=4.7truecm
\slika{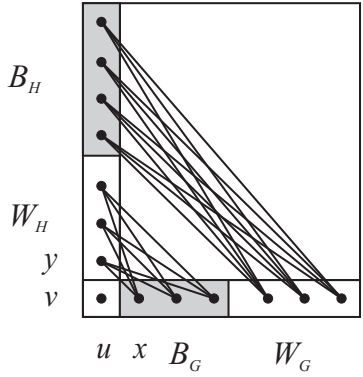}
\caption{The $S$-edges of $K(u,v)$.}
\label{okolice}
\end{figure}

We claim that $\alpha(u,v)\geq\delta(G)$ for all $(u,v)\in X$. Denote by 
$B_G,B_H,W_G$ and $W_H$ the size of black and white parts of $N_G(u)\times \{v\}$ 
and $\{u\}\times N_H(v)$ (see Fig.~\ref{okolice}). 
Since $(x,v)$ is black and $(u,y)$ is white (if $(u,v)\in W_G(x,y)\times B_H(x,y)$, otherwise the converse is true) 
we find that $B_G\geq 1$ and $W_H\geq 1$.  
Note also that $S$-edges of $K(u,v)$ incident to $(x,v)$ or $(u,y)$  
are not counted twice, and hence contribute to $\alpha''(u,v)$.
Therefore $$\alpha(u,v)\geq B_G+W_H-1+\frac12 ((B_G-1)(W_H-1)+B_HW_G)\,.$$ 
If $B_H=0$ or $W_G=0$, then $W_H\geq \delta(H)$ or 
$B_G\geq \delta(G)$,  and hence $\alpha(u,v)\geq \delta(H)$ or $\alpha(u,v)\geq \delta(G)$. 
So assume $B_H>0$ and $W_G> 0$. 
Since $B_G\geq \delta(G)-W_G$ and $W_H\geq \delta(H)-B_H$ we find that 
$$\alpha(u,v)\geq \frac12 (B_G+W_H)-1+\frac12(\delta(G)+\delta(H))+ \frac12 ((B_G-1)(W_H-1)+B_HW_G-B_H-W_G).$$
Since $B_H>0$ and $W_G> 0$ we see that $B_HW_G-B_H-W_G\geq -1$.
Therefore the only possibility for $\alpha(u,v)<\delta(G)$ is when 
$$B_HW_G-B_H-W_G=-1,\; B_G+W_H=2~{\rm and}~\delta(G)=\delta(H)\,.$$
If the above equalities are true we have $W_H=1$ and $B_G=1$. Since $W_H=1$ we find that 
$B_H\geq 2$ (for otherwise $\delta(H)\leq 2$) and therefore 
$B_HW_G-B_H-W_G=-1$ only if $B_H=2$ and $W_G=1$. But then $\delta(G)=2$ and hence 
$\delta(G)\neq \delta(H)$.
This proves that $\alpha(u,v)\geq\delta(G)$ for all $(u,v)\in X$, 
and therefore 
\begin{equation} \label{enacba}
|S|\geq \sum_{(u,v)\in X}\alpha(u,v)\geq \sum_{(u,v)\in X} \delta(G)\geq \delta(G)\delta(H) 
\end{equation}
and the equality holds only if $|X|=\delta(H)$ and $\alpha(u,v)=\delta(G)$ for all $(u,v)\in X$.

\begin{figure}[htb!]
\epsfxsize=12.0truecm
\slika{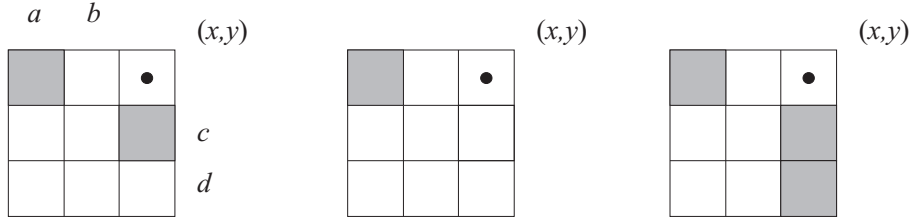}
\caption{When $\delta (H)=\delta (G)=2$.}
\label{delta2}
\end{figure}

Now assume that  $2\geq \delta (H)\geq \delta (G)$. If $\delta (H)=\delta (G)=2$ then 
$\alpha(u,v)\geq 3/2$ for every $(u,v)\in X$. So if $|X|\geq 3$ we have 
$$|S|\geq \sum_{(u,v)\in X}\alpha(u,v)\geq  9/2 >\delta(G)\delta(H) \,.$$ 
If $|X|=2$ we have one of the three cases presented in Fig. \ref{delta2}. Let $N(x)=\{a,b\}$ and 
$N(y)=\{c,d\}$. 

For the second and third case we have $\alpha(u,v)\geq 2$ for all $(u,v)\in X$ 
(because the set $X$ is a subset of a single $H$-layer and therefore $\alpha'(u,v)=0$), 
and thus $|S|\geq \sum_{(u,v)\in X}\alpha(u,v)\geq 4$. For the first case we argue that $\alpha(u,v)=2$ 
unless $(a,c)$ is adjacent to $(b,d)$ and one of these two vertices is black and the other white. 
Since $\alpha(b,c)+\alpha(a,d)\geq 3$ and there is one additional $S$-edge (either 
$(x,y)(a,c)$ or $(x,y)(b,d)$) we find that $|S|\geq 4$. Finally if  $\delta(G)=1$ we have 
$\alpha (u,v)\geq \delta (G)$ and since $|X|\geq \delta(H)$ we conclude  $|S|\geq \delta(G)\delta(H)$.

\qed

Several details of the above proof are needed in the following section where we describe
the structure of a minimum edge-cut in a direct product of graphs. We use the notation of the above proof and 
write these details in remarks that follow. Additionaly we define  
$$T=S\cap \{e\in E(K(u,v))\,|\,(u,v)\in X\}$$ 
and use this notation in the rest of the paper.

\begin{remark} \label{edina2}
Let $G$ and $H$ be graphs with $\delta(G)\leq \delta (H)$ and $\delta(H)\geq 3$. Suppose that 
$(u,v)\in X, \alpha'(u,v)=0$ and $\alpha(u,v)=\delta(G)$. 
\begin{itemize}
\item[(i)] If $(u,v)\in W_G(x,y)\times B_H(x,y)$ then $W_H=1$ and $W_G=0$, or $B_G=1$ and $B_H=0$. 
Moreover if $\delta(G)<\delta(H)$ then $W_H=1$ and $W_G=0$.
\item[(ii)] If $(u,v)\in B_G(x,y)\times W_H(x,y)$ then $B_H=1$ and $B_G=0$, or $W_G=1$ and $W_H=0$. 
Moreover if $\delta(G)<\delta(H)$ then $B_H=1$ and $B_G=0$.
\end{itemize}
\end{remark}

\proof Since  $\alpha'(u,v)=0$ we find that $\alpha(u,v)=B_GW_H+W_GB_H$. 
The result follows from $\delta(G)\leq \delta(H)$ and $\delta(H)\geq 3$.
\qed

The following remark follows directly from inequality \eqref{enacba} and the fact that 
$|X|\geq \delta(H)$ and $\alpha(u,v)\geq \delta(G)$.

\begin{remark}\label{stefka}
Let $G$ and $H$ be graphs with $\delta(G)\leq \delta (H)$ and $\delta(H)\geq 3$, and let $S$ be a minimum edge-cut in 
$G\times H$. If there is a vertex  $(x,y)\in V(G\times H)$ such that $W_G(x,y)\neq \emptyset$ and $B_G(x,y)\neq \emptyset$ then 
$|X|=\delta(H)$ and $\alpha(u,v)=\delta(G)$ for every $(u,v)\in X$.
\end{remark}

In  proof of Lemma \ref{glavna} we have showed that $|T|\geq \delta(G)\delta(H)$ except in the 
case when $\delta(G)=\delta(H)=2$ where we needed an additional $S$-edge to prove that 
$|S|\geq \delta(G)\delta(H)$. We write this in the following remark.
\begin{remark} \label{jozica}
Let $G$ and $H$ be graphs such that $\delta(G)\leq \delta (H)$ and $\delta(H)\geq 3$. 
If there is a vertex $(x,y)\in V(G\times H)$ such that 
$B_G(x,y)\neq \emptyset$ and $W_G(x,y)\neq \emptyset$ 
then $|T|\geq \delta(G)\delta(H)$.

\end{remark}

For a graph $G$ and $A,B\subseteq V(G)$, 
we denote by $[A,B]$ the set of all
vertices with one endvertex in $A$ and the other in $B$. Additionaly, 
we denote by $[A]$ the subgraph of $G$ induced by $A$. 

Let $G$ and $H$ be graphs. 
If $S$ is a minimum edge cut of $G$ and $X,Y$ are connected components of $G-S$  
then $$[X\times V(H),Y\times V(H)]$$
is an edge cut of $G\times H$ whose connected components are 
$X\times V(H)$ and $Y\times V(H)$ as shown in case 1 of
Fig.~\ref{struktura}. Clearly  
$$|[X\times V(H),Y\times V(H)]|=2\lambda(G)|E(H)|\,,$$
and hence $\lambda(G\times H)\leq 2\lambda(G)|E(H)|.$ Similarely we prove 
$\lambda(G\times H)\leq 2\lambda(H)|E(G)|$, and obviously  
$\lambda(G\times H)\leq \delta(G\times H)$.

Recall that the bipartite edge frustration  of a graph $G$, denoted as $\varphi (G)$, is the 
minimum number of edges whose delition makes $G$ a bipartite graph. 
Denote by $\mathcal{P}(V(G))$ the set of all partitions of $V(G)$ on two  subsets 
(one of them may be empty), 
and let 
$$\rho (G)=\min \{ 2\varphi([A])+|[A,B]|:\{A,B\}\in \mathcal{P}(V(G)), A\neq \emptyset\}\,.$$
For graphs $G$ and $H$ we define 
$$\psi(G,H)=\rho(G)|E(H)|+2\varphi(H)|E(G)|\,.$$ 

\begin{figure}[htb!]
\epsfxsize=4.6truecm
\slika{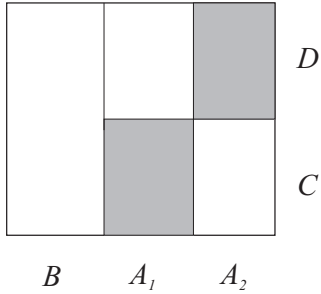}
\caption{$\lambda(G\times H)\leq \psi(G,H)$.}
\label{psi}
\end{figure}

We claim that $\lambda(G\times H)\leq\psi(G,H)$. 
To see this let $\{A,B\}\in \mathcal{P}(V(G))$ be a partition of $V(G)$
for which $2\varphi([A])+|[A,B]|$ is minimum, 
and let $C\cup D$ be a partition of 
$V(H)$ such that  $|[C,D]|=|E(H)|-\varphi(H)$. 
Assume also that $A=A_1\cup A_2$ is a partition of $A$ with $|[A_1,A_2]|=|E([A])|-\varphi([A])$ 
(so there are $\varphi([A])$ edges with both endvertices in $A_1$ or both in $A_2$). Let 
$(A_1\times C)\cup (A_2\times D)$ be black, and $(A_2\times C)\cup (A_1\times D)\cup (B\times V(H))$ 
 white,  see Fig. \ref{psi}.  We will show that the number of edges between 
black and white part of $V(G\times H)$ is at most 
$\rho(G)|E(H)|+2\varphi(H)|E(G)|$.

The number of edges with both endvertices in $V(G)\times C$ or 
both in $V(G)\times D$ is at most $2\varphi(H)|E(G)|$. So we only need to count the 
edges with one endvertex in $V(G)\times C$ and the other in $V(G)\times D$. 
For each edge in $[A,B]$ and each edge in $[C,D]$ there is exactly one edge in $G\times H$ 
with one endvertex black and the other  white, so the total number of such edges is 
at most $|E(H)||[A,B]|$. For each edge with both endvertices in $A_1$ (or both in $A_2$), 
and each edge in $[C,D]$, there are two edges in $G\times H$ 
with one endvertex black and the other  white. Therefore the total number of such edges is at most 
$2\varphi([A])|E(H)|$. This proves that $\lambda(G\times H)\leq\psi(G,H)$. If $H$ is bipartite then 
$[C,D]=V(H)$ so the number of $S$-edges with one endvertex in $V(G)\times C$ and the other in $V(G)\times D$ 
is equal to $\rho(G)|E(H)|$. This explains the following remark. 

\begin{remark} \label{45}
If $H$ is bipartite with $E(H)=[C,D]$, then $\psi(G,H)$ is equal to the minimum size of an 
edge-cut of type 4 or 5 such that $(A_1\times C)\cup(A_2\times D)$ is black (and the rest white) 
for some $A_1,A_2\subseteq V(G)$. 
Analagously, if $G$ is bipartite with $E(G)=[A,B]$, then $\psi(H,G)$ 
is equal to the minimum size of an edge-cut of type 4 or 6 such that $(A\times C_1)\cup (B\times C_2)$ is black 
(and the rest white) for some $C_1,C_2\subseteq V(H)$. 
\end{remark}

Combining all upper bounds into one we get  
$$\lambda (G\times H)\leq \min \{2\lambda(G)|E(H)|,2\lambda(H)|E(G)|,\delta(G\times H), 
\psi(G,H), \psi(H,G)\}\,.$$ 
We prove next that the above ineqality is in fact an equality. \\

\noindent {\bf Proof of Theorem \ref{formula}.}
If $G$ or $H$ is not connected then both sides are 0. So assume that 
$G$ and $H$ are connected. 
Let $S$ be a minimum edge-cut in $G\times H$ and $B,W$ be connected components of 
$(G\times H)-S$. We need to prove that 
$$|S|\geq \min \{2\lambda(G)|E(H)|,2\lambda(H)|E(G)|,\delta(G\times H), \psi(G,H), \psi(H,G)\}\,.$$ 

If $|S|\geq \delta(G\times H)$ then we are done. So assume that $|S|< \delta(G\times H)$.   
Without loss of generality assume also $\delta(G)\leq \delta(H)$. By Lemma \ref{glavna}  
for every $(x,y)\in V(G\times H)$ either 
$N(x)\times \{y\}\subseteq B$ or $N(x)\times \{y\}\subseteq W$. 
Let $y\in V(H)$ be any  vertex. If there are two 
adjacent vertices $x_1,x_2\in V(G)$ such that $(x_1,y)$ and $(x_2,y)$ are both black (or both white) then  
we find that $G_y\subseteq B$ (or $G_y\subseteq W$). 
If there are no such vertices, then we have a proper 2-coloring of $G$, and hence $G$ is bipartite.  Moreover 
if $X\cup Y$ is the unique (recall that $G$ is connected) bipartition of $V(G)$ 
then either $X\times \{y\}\subseteq W$ and $Y\times \{y\}\subseteq B$, or 
$X\times \{y\}\subseteq B$ and $Y\times \{y\}\subseteq W$. Since this is true for every $y\in V(H)$ 
wee find that the edge-cut $S$ is of type 2,4,6, or 8 (see Fig. \ref{struktura}). 
Moreover types 4,6 and 8 can occure only if $G$ is bipartite. 
If $S$ is of type 2 we see that $|S|\geq 2\lambda(H)|E(G)|$. Otherwise $G$ is bipartite 
and therefore $\varphi(G)=0$. It follows that $\psi(H,G)=\rho(H)|E(G)|$. 
In a partition of type 8 we may change the color of all black $G$-layers to become white. 
In this way we get a partition of type  6, moreover such partition 
requires deletion of fewer or equal number of edges (because $G$ is bipartite). 
If $S$ is of type 4 or 6 we have $|S|\geq \rho(H)|E(G)|=\psi(H,G)$ (see Remark \ref{45}), 
and therefore the same applies for type 8 (because it requires deletion of at least as much edges as type 6). 
\qed

The proof of the above theorem gives  a more precise information regarding the structure 
of connected components. We write this in the following corollary.

\begin{corollary}\label{cucek}
Let $G$ and $H$ be graphs with $\delta(G)\leq \delta(H)$. 
If $S$ is a minimum edge-cut in $G\times H$ and 
 for every $(x,y)\in V(G\times H)$ either 
$N(x)\times \{y\}\subseteq B$ or $N(x)\times \{y\}\subseteq W$
then $S$ is of type 2,4,6 or 8.
\end{corollary}

\section{The structure of a minimum edge cut}
\label{zadnje}
In this section we use the notation and definitions of the previous section. \\

\noindent {\bf Proof of Theorem \ref{tipi}.}
Let $S$ be a minimum edge-cut in $G\times H$ and $B,W$ be connected components of $(G\times H)-S$.
Assume that $\delta(G)\leq \delta(H)$ and $\delta(H)\geq 3$.

Corollary \ref{cucek} settles the structure  of a minimum edge-cut if 
$N(x)\times \{y\}\subseteq B$ or $N(x)\times \{y\}\subseteq W$ for every $(x,y)\in V(G\times H)$. 
Assume now that there is $(x,y)\in V(G\times H)$ such that 
$$B_G(x,y)=(N(x)\times \{y\})\cap B\neq \emptyset ~{\rm and}~ W_G(x,y)=(N(x)\times \{y\})\cap W\neq \emptyset\,.$$

\begin{figure}[htb!]
\epsfxsize=15.0truecm
\slika{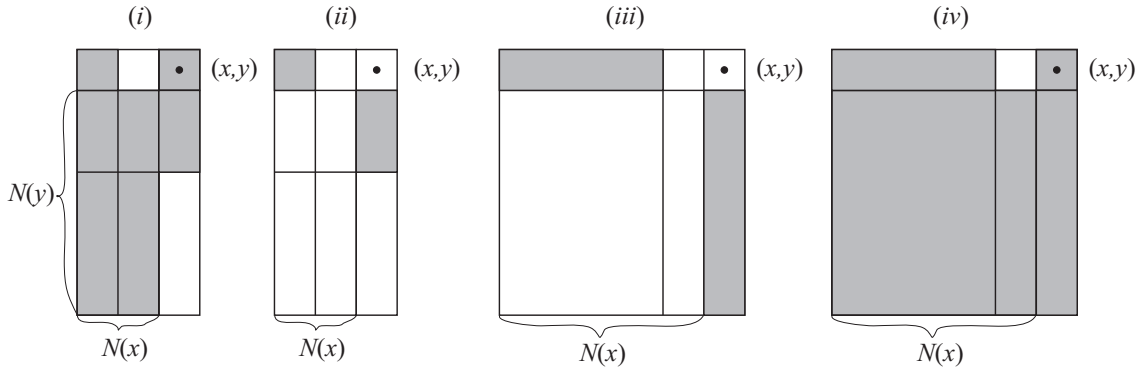}
\caption{The four cases when $B_G(x,y)\neq \emptyset$ and $W_G(x,y)\neq \emptyset\,$.}
\label{macka}
\end{figure}

By Remark \ref{stefka} we know that 
$|X|=\delta(H)$ and $\alpha(u,v)=\delta(G)$ for every $(u,v)\in X$. It follows from 
Remark \ref{jozica} that
$|T|\geq \delta(G)\delta(H)$ 
and therefore $S=T$. 
This is particularly noteworthy since most of the arguments in the proof below reduce 
to the fact that there are no $S$-edges other than the edges in $T$.

Since $|X|=\delta(H)$ we see that one of the four cases 
shown in Fig. \ref{macka} occures.
Note that in all four cases the color of vertex $(x,y)$ is the same as the color of 
vertices in $N_{G\times H}(x,y)$ because no edge incident to $(x,y)$ is in 
$T$ (and therefore also not in $S$). 
Therefore for every $(u,v)\in X$ we have $\alpha'(u,v)=0$. Since also $\alpha(u,v)=\delta(G)$ for 
every $(u,v)\in X$ we may use Remark \ref{edina2} in the sequel. Observe also that 
\begin{equation}\label{min}
\alpha(u,v)\geq \min\{\deg_G(u),\deg_H(v)\}\,.
\end{equation}

\begin{figure}[htb!]
\epsfxsize=5.5truecm
\slika{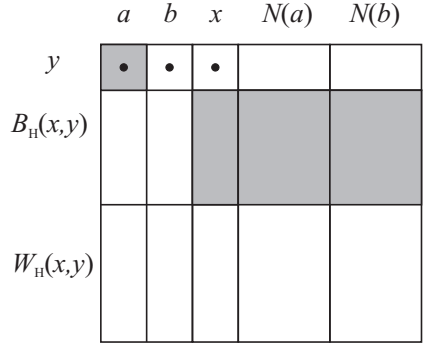}
\caption{Cases $(i)$ and $(ii)$.}
\label{micka}
\end{figure}

{\em Case (i) and (ii).} The color of vertex $(x,y)$ is irrelevant, so we prove both cases simultaneously.
We claim that $G$ is $P_3$ or $C_4$. Let $N(x)=\{a,b\}$. Consider open neighborhoods of vertices $a$ and 
$b$ as shown in Fig. \ref{micka}, and note that $x\in N(a)\cap N(b)$. We claim that 
$(N(a)\cup N(b))\times B_H(x,y)$ is  black and $(N(a)\cup N(b))\times W_H(x,y)$ 
is  white. Since every vertex in 
$N(a)\times B_H(x,y)$ is adjacent to $(a,y)$, and the corresponding edges are not in   
$T$, we find that $N(a)\times B_H(x,y)$ is black. 
Since $\delta(G)\leq 2$ we find that $\delta(H)>\delta(G)$ and therefore
$N(b)\times B_H(x,y)$ is black by Remark~\ref{edina2} (because $W_H=1$ and $W_G=0$ for 
every $(u,v)\in W_G(x,y)\times B_H(x,y)$).
Similar arguments prove that $(N(a)\cup N(b))\times W_H(x,y)$ is white. 

Next we claim that $N(N(a)\cup N(b))\subseteq \{a,b\}$. If not, there is a   
$t\neq a,b$ adjacent to a vertex $t'\in N(a)\cup N(b)$. Hence  
$(t,y)$ is adjacent to  black and white vertices in $\{t'\}\times B_H(x,y)$ and 
$\{t'\}\times W_H(x,y)$, respectively. Thus there is at least one $S$-edge
not in $T$, a contradicition. 
If $N(a)=N(b)=\{x\}$ then $G=P_3$. If $N(a)\neq \{x\}$ 
or $N(b)\neq \{x\}$ then $N(a)=N(b)$, for otherwise $\delta(G)=1$ and hence $\alpha(u,v)>\delta(G)$ for 
a $(u,v)\in X$ (see \eqref{min}). Since $\delta(G)\leq 2$ we find that $|N(a)|\leq 2$ 
(for otherwise inequality \eqref{min} implies $\alpha(u,v)\geq 3>\delta(G)$ for all $(u,v)\in X$).
Since $a$ and $b$ are nonadjacent and the vertices in 
$N(a)$ are nonadjacent (otherwise we have $S$-edges that are not in $T$) we find that $G$ is either $P_3$ or $C_4$.

{\em Case (iii)}.
Let $u\in B_G(x,y)$ and observe that $H_x$ and $H_u$ are adjacent and in both layers there are
black and white vertices. 
Since $S=T$ and hence $[H_x,H_u]\cap S=\emptyset$ we find that the subgraph of $G\times H$ induced by $H_x\cup H_u$ 
is not connected, and therefore $H$ is bipartite by Observation \ref{kuku}. 
Let $W_G(x,y)=\{b\}$ and $E(H)=[C,D]$.
Consider the graph $G-b$ (vertex deleted subgraph of $G$), 
and observe that no $S$-edge  has both endvertices in $V(G-b)\times V(H)$ (because $S=T$ and each $T$-edge has one endvertex in $H_b$).  
Suppose that $G-b$ has no isolated vertices 
and consider two adjacent $H$-layers $H_s$ and $H_t$ where $s,t\neq b$. 
By  Observation \ref{kuku} one of the following occures:

\begin{quote}
$H_s$ and $H_t$ are both white \\
$H_s$ and $H_t$ are both black \\
$(\{s\}\times C)\cup (\{t\}\times D)$ is white and $(\{s\}\times D)\cup (\{t\}\times C)$ is black \\
$(\{s\}\times C)\cup (\{t\}\times D)$ is black and $(\{s\}\times D)\cup (\{t\}\times C)$ is white 
\end{quote}
The layers $H_u$ and $H_x$ fall into the last two cases.
Therefore each $H$-layer, except $H_b$, is either black, or white, or partitioned 
(into black and white parts) like $H_x$, or like $H_u$. Now the $H_b$ layer may only be adjacent to layers 
that are white or partitioned like $H_x$ (because if $b'$ is adjacent to $b$, then $(b',y)$ is white 
unless there is an $S$-edge not in $T$). 
By Remark \ref{edina2} either all layers adjacent to $H_b$ are white or all are partitioned as $H_x$. 
Moreover there are no black $H$-layers. 
In the first case $H_b$ is white and $S$ is of type 5, and in the second case 
$G$ is bipartite (since also $H$ is bipartite this is a contradiction). Finally if $G-b$ has an isolated vertex then 
$\delta(G)=1$ and therefore $N(b)=\{x\}$ by \eqref{min}. In this case $G$ is bipartite, a contradiction.

\begin{figure}[b!]
\epsfxsize=13.0truecm
\slika{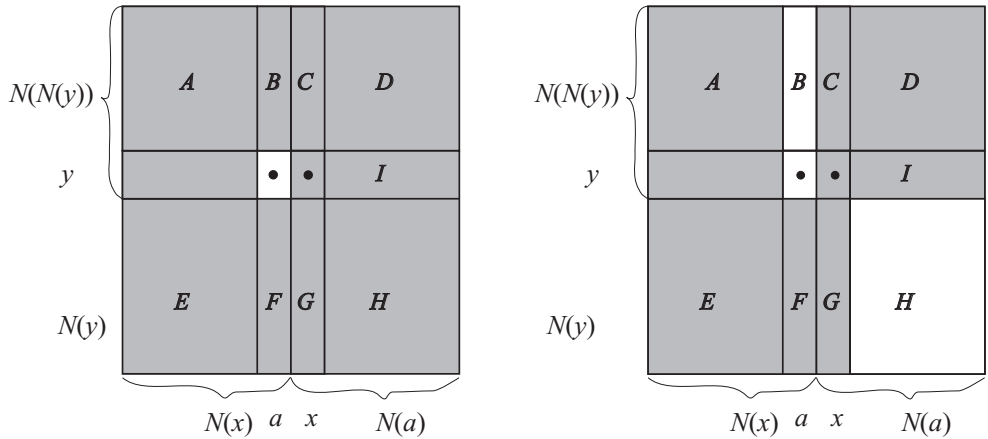}
\caption{Possibilities in case $(iv)$.}
\label{mjav}
\end{figure}

{\em Case (iv)}. Here we have two cases (see Fig.~\ref{mjav}). In any case, the vertex $(x,y)$ is black, and 
so are also parts $E$ and $F$. Since $A$ is adjacent to $G$, and $G$ is black, we find that also 
$A$ is black. Parts $C,D$ and $I$ are black because they are adjacent to $F$. The only difference 
of these two cases is in parts $B$ and $H$. Either both are black or both white (see Remark \ref{edina2}). 
If both are black, then there is only one white vertex 
and hence $|W|=1$, otherwise $S$ is of type 5. This proves the theorem if $\delta(H)\geq 3$.

Now we consider the case $\delta(G)\leq \delta(H)\leq 2$.  
We prove it by a case analysis. 
Assume that 
$\delta(G)=\delta(H)=2$. If $|X|\geq 3$ then 
$$|S|\geq \sum_{(u,v)\in X} \alpha(u,v)\geq 3\cdot \frac 32>4=\delta(G)\delta(H)\,,$$
because $\alpha(u,v)\geq 3/2$ for all $(u,v) \in X.$ So assume that $|X|=2$. 
One of the cases shown in Fig. \ref{delta2} occures. 
Let $N(x)=\{a,b\}$ and $N(y)=\{c,d\}$. Since 
$\delta(G)=\delta(H)=2$ there is at least one neighbor $a'\neq x$ of $a$, 
$b'\neq x$ of $b$, $c'\neq y$ of $c$, and $d'\neq y$ of $d$.  

Consider the first case where 
$(a,y)$ and $(x,c)$ are black, and $(b,y)$ and $(x,d)$ are white. 
Since $(a',c)$ is adjacent to $(a,y)$ we find that $(a',c)$ is black. 
If $(a',d)$ is  white then 
$N(a')\subseteq \{a,b\}$, and since $\delta(G)=2$ we find that $N(a')=\{a,b\}$. If  
$a'\neq b'$ then 
$\{x,a',b'\}\subseteq N(b)$ and hence $\alpha(b,c)\geq 3>\delta(G)$, a contradiction. Otherwise, 
if $a'=b'$, then $G=C_4$.
Assume now that both $(a',c)$ and $(a',d)$ are black (note that in this case $a'\neq b'$ because of the edge $(b',d)(b,y)$). 
If $N(a')\subseteq \{a,b\}$ then we argue the same as before.
If not,  there is a vertex $a''\neq a,b$ adjacent to $a'$. 
By Observation \ref{kuku} $H$ is bipartite (because $a'$-layer has black and white vertices, and no $T$-edge 
has one endvertex in $a'$-layer and the other in $a''$-layer). Similar arguments prove that 
either $\alpha(a,d)\geq 3$ (if $c'$ is adjacent to $d$) or $G$ is biparite. 
In both cases we have a contradiction (either both graphs are bipartite or 
$\alpha(u,v)>\delta(G)$ for a $(u,v)\in X$).

Consider now the second case  of Fig. \ref{delta2}. 
We claim that $S$ is of type 3 or $H$ is $C_4$. Since $\delta(G)=2$ we find that $a'$ has a neighbor $a''\neq a$.
If $(a',c)$ is black and $(a',d)$ is white, or vice versa, we have an $S$-edge incident to $(a'',y)$, a contradiction.
Therefore either $(a',c)$ and $(a',d)$ are both black or both white. 
If they are both black, then $H$ is bipartite.
Since $G$ is not bipartite we find that layers $G_{c'}$ and $G_{d'}$ 
are both white unless $c'=d'$. If $c'=d'$ then $H=C_4$, otherwise $(a',c')$ and $(a',d')$ are both white, 
which is in a contradiction to Observation \ref{kuku}. 
If $(a',c)$ and $(a',d)$ are both white then $S$ is of type 3. 

In third case of Fig. \ref{delta2} the graph $H$ is bipartite. 
If $H\neq C_4$ we find that $c'$ has a neighbor $c''\neq c$. 
Since  $G$ is not bipartite we find that layers $G_{c'}$ and $G_{c''}$ are both white and hence 
$(x,c')$ and $(x,c'')$ are white, a contradiction.

Assume that $\delta(G)=1$ and $\delta(H)=2$. If 
$(x,c)$ is black  and  $(x,d)$ white we find that $G=P_3$ (because $N(x)=\{a,b\}$). 
If $(x,c)$ and  $(x,d)$ are both black then $H$ is bipartite. If $c'=d'$ then $H=C_4$, otherwise also $G$ 
is bipartite, a contradiction. 

Finally we look at the case $\delta(G)=\delta(H)=1$. 
Assume that $|V(G)|,|V(H)|\geq 3$.  If all but one $G$-layer (or $H$-layer) are  white 
we find that $S$ is of type 3. Otherwise both graphs are bipartite. 
If $H=K_2$ we find that $S$ is of type 3 or 5.

\qed

\end{document}